\documentclass{amsart}
       \usepackage{amsmath,amsthm,upref}
       \usepackage{amssymb,amsbsy}
       \usepackage{amsfonts}
\usepackage{times}
\theoremstyle{definition}

\newtheorem{definition}{Definition}
\newtheorem{corollary}[definition]{Corollary}

\theoremstyle{plain}
\newtheorem{theorem}[definition]{Theorem}
\newtheorem{proposition}[definition]{Proposition}

\newtheorem{lemma}[definition]{Lemma}
\theoremstyle{remark}
\newtheorem*{remark}{REMARK}
\newtheorem*{example}{Example}

\begin{document}

\author{Andrzej Borowiec}
\address{Institute of Theoretical Physics, University of Wroc{\l }aw, Pl. Maxa Borna 9,
50-204 Wroc{\l }aw, Poland}
\email{borow@ift.univ.wroc.pl}
\author{Wies{\l }aw A. Dudek}
\address{Institute of Mathematics, Technical University of Wroc{\l }aw, Wybrzeze
Wyspianskiego 27, 50-370 Wroc{\l }aw, Poland}
\email{dudek@im.pwr.wroc.pl}
\author{Steven Duplij}
\address{Department of Physics and Technology, Kharkov National University, Kharkov
61001, Ukraine}
\email{Steven.A.Duplij@univer.kharkov.ua}
\urladdr{http://www.math.uni-mannheim.de/\symbol{126}duplij}

\title{Basic concepts of ternary Hopf algebras}
\date{August 29, 2001\\
\mbox{}\, \, \, {\it Published in}: Journal of Kharkov National University, ser. Nuclei, Particles and Fields, 
{\bf 529}, N 3(15) (2001) pp. 21--29. Reprints available from SD}
\maketitle
\begin{abstract}
\noindent
The theory of ternary semigroups, groups and algebras is reformulated in the
abstract arrow language. Then using the reversing arrow ansatz we define
ternary comultiplication, bialgebras and Hopf algebras and investigate their properties.
The main property "to be binary derived" is considered in detail. The co-analog of Post theorem
is formulated.
It is shown that there exist 3 types of ternary coassociativity, 3 types of ternary counits
and 2 types of ternary antipodes. Some examples are also presented.

\end{abstract}

Ternary and $n$-ary generalizations of algebraic structures is the most
natural way for further development and deeper understanding of their
fundamental properties. Firstly ternary algebraic operations were introduced
already in the XIX-th century by A. Cayley. As the development of Cayley's
ideas it were considered $n$-ary generalization of matrices and their
determinants \cite{sokolov,kap/gel/zel} and general theory of $n$-ary algebras
\cite{law,car5} and ternary rings \cite{lis} (for physical applications in
Nambu mechanics, supersymmetry, Yang-Baxter equation, etc. see
\cite{ker1,vai/ker} as surveys). The notion of an $n$-ary group was introduced
in 1928 by W. D\"{o}rnte \cite{dor3} (inspired by E. N\"{o}ther) which is a
natural generalization of the notion of a group and a ternary group considered
by Certaine \cite{cer} and Kasner \cite{kas}.
%
%
For many applications of $n$-ary groups and quasigroups see
\cite{rusakov1,cup/cel/mar/dim} and \cite{belousov} respectively. From another
side, Hopf algebras \cite{abe,sweedler} and their generalizations
\cite{nil1,nik/vai,dup/li1,dup/li2} play a basic role in the quantum group
theory (see e.g. \cite{demidov,kassel,shn/ste}).

In the first part of this paper we reformulate necessary material on ternary
semigroups, groups and algebras \cite{belousov,rusakov1} in the abstract arrow
language. Then according to the general scheme \cite{abe} using systematic
reversing order of arrows, we define ternary bialgebras and Hopf algebras,
investigate their properties and present examples.

\section*{Ternary semigroups}

A non-empty set $G$ with one \emph{ternary} operation $[\;]:G\times G\times
G\rightarrow G$ is called a \emph{ternary groupoid} and is denoted by
$(G,[\;])$ or $\left(  G,m^{\left(  3\right)  }\right)  $. We will present
some results using second notation, because it allows to reverse arrows in the
most clear way. In proofs we will mostly use the first notation due to
convenience and for short.

If on $G$ there is a binary operation $\odot$ (or $m^{\left(  2\right)  }$)
such that $[xyz]=(x\odot y)\odot z$ or%
\begin{equation}
m^{\left(  3\right)  }=m_{der}^{\left(  3\right)  }=m^{\left(  2\right)
}\circ\left(  m^{\left(  2\right)  }\times\operatorname*{id}\right)
\label{mder}%
\end{equation}
for all $x,y,z\in G$, then we say that $[\;]$ or $m_{der}^{\left(  3\right)
}$ is \emph{derived} from $\odot$ or $m^{\left(  2\right)  }$and denote this
fact by $(G,[\;])=der(G,\odot)$. If
\[
\lbrack xyz]=((x\odot y)\odot z)\odot b
\]
holds for all $x,y,z\in G$ and some fixed $b\in G$, then a groupoid $(G,[\;]$
is \emph{$b$-derived} from $(G,\odot)$. In this case we write
$(G,[\;])=der_{b}(G,\odot)$ (cf. \cite{dud/mic1,dud/mic2}).
%
%

We say that $(G,[\;]$ is a \emph{ternary semigroup} if the operation $[\;]$ is
\emph{associative}, i.e. if
\begin{equation}
\left[  \left[  xyz\right]  uv\right]  =\left[  x\left[  yzu\right]  v\right]
=\left[  xy\left[  zuv\right]  \right]  \label{ass1}%
\end{equation}
holds for all $x,y,z,u,v\in G$, or%
\begin{equation}
m^{\left(  3\right)  }\circ\left(  m^{\left(  3\right)  }\times
\operatorname*{id}\times\operatorname*{id}\right)  =m^{\left(  3\right)
}\circ\left(  \operatorname*{id}\times m^{\left(  3\right)  }\times
\operatorname*{id}\right)  =m^{\left(  3\right)  }\circ\left(
\operatorname*{id}\times\operatorname*{id}\times m^{\left(  3\right)
}\right)  \label{assm}%
\end{equation}

Obviously, a ternary operation $m_{der}^{\left(  3\right)  }$ derived from a
binary associative operation $m^{\left(  2\right)  }$ is also associative in
the above sense, but a ternary groupoid $(G,[\;])$ $b$-derived ($b$ is a
cansellative element) from a semigroup $(G,\odot)$ is a ternary semigroup if
and only if $b$ lies in the center of $(G,\odot)$.

Fixing in a ternary operation $m^{\left(  3\right)  }$ one element $a$ we
obtain a binary operation $m_{a}^{\left(  2\right)  }$. A binary groupoid
$(G,\odot)$ or $\left(  G,m_{a}^{\left(  2\right)  }\right)  $, where $x\odot
y=[xay]$ or%
\begin{equation}
m_{a}^{\left(  2\right)  }=m^{\left(  3\right)  }\circ\left(
\operatorname*{id}\times a\times\operatorname*{id}\right)  \label{ma}%
\end{equation}
for some fixed $a\in G$ is called a \emph{retract} of $(G,[\;])$ and is
denoted by $ret_{a}(G,[\;])$. In some special cases described in
\cite{dud/mic1,dud/mic2} we have $(G,\odot)=ret_{a}(der_{b}(G,\odot))$ or
$(G,\odot)=der_{c}(ret_{d}(G,[\;]))$.

\begin{lemma}
If in the ternary semigroup $(G,[\;])$ or $\left(  G,m^{\left(  3\right)
}\right)  $ there exists an element $e$ such that for all $y\in G$ we have
$\left[  eye\right]  =y$, then this semigroup is derived from the binary
semigroup $\left(  G,m_{e}^{\left(  2\right)  }\right)  $, where%
\begin{equation}
m_{e}^{\left(  2\right)  }=m^{\left(  3\right)  }\circ\left(
\operatorname*{id}\times e\times\operatorname*{id}\right)  \label{me}%
\end{equation}
In this case $(G,[\;])=der(ret_{e}(G,[\;])$.
\end{lemma}

\begin{proof}
Indeed, if we put $x\circledast y=[xey]$, then $(x\circledast y)\circledast
z=[[xey]ez]=[x[eye]z]=[xyz]$ and $x\circledast(y\circledast
z)=[xe[yez]]=[x[eye]z]=[xyz]$, which completes the proof.
\end{proof}

The same ternary semigroup $\left(  G,m^{\left(  3\right)  }\right)  $ can be
derived from two different semigroups $\left(  G,\circledast\right)  $ or
$\left(  G,m_{e}^{\left(  2\right)  }\right)  $ and $\left(  G,\diamond
\right)  $ or $\left(  G,m_{a}^{\left(  2\right)  }\right)  $. Indeed, if in
$G$ there exists $a\neq e$ such that $[aya]=y$ for all $y\in G$, then by the
same argumentation we obtain $[xyz]=x\diamond y\diamond z$ for $x\diamond
y=[xay]$. In this case for $\varphi(x)=x\diamond e=[xae]$ we have
\[
x\circledast y=[xey]=[x[aea]y]=[[xae]ay]=(x\diamond e)\diamond y=\varphi
(x)\diamond y
\]
and
\[
\varphi(x\circledast y)=[[xey]ae]=[[x[aea]y]ae]=[[xae]a[yae]]=\varphi
(x)\diamond\varphi(y).
\]
Thus $\varphi$ is a binary homomorphism such that $\varphi(e)=a$. Moreover for
$\psi(x)=[eax]$ we have
\begin{align*}
\psi(\varphi(x))  &  =[ea[xae]]=[e[axa]e]=x,\\
\varphi(\psi(x))  &  =[[eax]ae]=[e[axa]e]=x
\end{align*}
and
\[
\psi(x\diamond y)=[ea[xay]]=[ea[x[eae]y]]=[[eax]e[aey]]=\psi(x)\circledast
\psi(y).
\]
Hence semigroups $(G,\circledast)$ and $(G,\diamond)$ are isomorphic.

\begin{definition}
An element $e\in G$ is called a \emph{middle identity} or a \emph{middle
neutral element} of $(G,[\;])$ if for all $x\in G$ we have $[exe]=x$ or%
\begin{equation}
m^{\left(  3\right)  }\circ\left(  e\times\operatorname*{id}\times e\right)
=\operatorname*{id}. \label{mmx}%
\end{equation}
An element $e\in G$ satisfying the identity $[eex]=x$ or%
\begin{equation}
m^{\left(  3\right)  }\circ\left(  e\times e\times\operatorname*{id}\right)
=\operatorname*{id}. \label{ml}%
\end{equation}
is called a \emph{left identity} or a \emph{left neutral element} of
$(G,[\;])$. By the symmetry we define a \emph{right identity}. An element
which is a left, middle and right identity is called a \emph{ternary identity}
(briefly: identity).
\end{definition}

There are ternary semigroups without left (middle, right) neutral elements,
but there are also ternary semigroups in which all elements are identities
\cite{rusakov1,pos}.

\begin{example}
In ternary semigroups derived from the symmetric group $S_{3}$ all elements of
order 2 are left and right (but no middle) identities.
\end{example}

\begin{example}
In ternary semigroup derived from Boolean group all elements are ternary
identities, but ternary semigroup $1$-derived from the additive group
$\mathbb{Z}_{4}$ has no left (right, middle) identities.
\end{example}

\begin{lemma}
For any ternary semigroup $(G,[\;])$ with a left (right) identity there exists
a binary semigroup $(G,\odot)$ and its endomorphism $\mu$ such that
\[
\lbrack xyz]=x\odot\mu(y)\odot z
\]
for all $x,y,z\in G$.
\end{lemma}

\begin{proof}
Let $e$ be a left identity of $(G,[\;])$. It is not difficult to see that the
operation $x\odot y=[xey]$ is associative. Moreover, for $\mu(x)=[exe]$, we
have
\[
\mu(x)\odot\mu(y)=[[exe]e[eye]]=[[exe][eey]e]=[e[xey]e]=\mu(x\odot y)
\]
and
\[
\lbrack xyz]=[x[eey][eez]]=[[xe[eye]]ez]=x\odot\mu(y)\odot z.
\]

The case of right identity the proof is analogous.
\end{proof}

\begin{definition}
We say that a ternary groupoid $(G, [\; ])$ is:

a \emph{left cancellative} if $[abx] = [aby]\Longrightarrow x=y$,

a \emph{middle cancellative} if $[axb] = [ayb]\Longrightarrow x=y$,

a \emph{right cancellative} if $[xab] = [yab]\Longrightarrow x=y$%
\newline holds for all $a,b\in G$.

A ternary groupoid which is left, middle and right cancellative is called
\emph{cancellative}.
\end{definition}

\begin{theorem}
\label{theor-canc}A ternary groupoid is cancellative if and only if it is a
middle cancellative, or equivalently, if and only if it is a left and right cancellative.
\end{theorem}

\begin{proof}
Assume that a ternary semigroup $(G,[\;])$ is a middle cancellative and
$[xab]=[yab]$. Then $[ab[xab]]=[ab[yab]]$ and in the consequence
$[a[bxa]b]=[a[bya]b]$ which implies $x=y$.

Conversely if $(G,[\;])$ is a left and right cancellative and $[axb]=[ayb]$
then $[a[axb]b]=[a[ayb]b]$ and $[[aax]bb]=[[aay]bb]$ which gives $x=y$.
\end{proof}

The above theorem is a consequence of the general result proved in
\cite{dud1}.
%

\begin{definition}
A ternary groupoid $(G,[\;])$ is \emph{semicommutative} if $[xyz]=[zyx]$ for
all $x,y,z\in G$. If the value of $[xyz]$ is independent on the permutation of
elements $x,y,z$, viz.%
\begin{equation}
\left[  x_{1}x_{2}x_{3}\right]  =\left[  x_{\sigma\left(  1\right)  }%
x_{\sigma\left(  2\right)  }x_{\sigma\left(  3\right)  }\right]  \label{xxx}%
\end{equation}
or $m^{\left(  3\right)  }=m^{\left(  3\right)  }\circ\sigma$, then $(G,[\;])$
is a \emph{commutative} ternary groupoid. If $\sigma$ is fixed, then a ternary
groupoid satisfying (\ref{xxx}) is called $\sigma$-commutative.
\end{definition}

The group $S_{3}$ is generated by two transpositions; $(12)$ and $(23)$. This
means that $(G,[\;])$ is commutative if and only if $[xyz]=[yxz]=[xzy]$ holds
for all $x,y,z\in G$.

As a simple consequence of Theorem \ref{theor-canc} from \cite{dud2}
%
we obtain

\begin{corollary}
If in a ternary semigroup $(G,[\; ])$ satisfying the identity $[xyz]=[yxz]$
there are $a,b$ such that $[axb]=x$ for all $x\in G$, then $(G,[\; ])$ is commutative.
\end{corollary}

\begin{proof}
According to the above remark it is sufficient to prove that $[xyz]=[xzy]$. We
have
\[
\lbrack xyz]=[a[xyz]b]=[ax[yzb]]=[ax[zyb]]=[a[xzy]b]=[xzy].
\]
\end{proof}

Mediality in the binary case $\left(  x\odot y\right)  \odot\left(  z\odot
u\right)  =\left(  x\odot z\right)  \odot\left(  y\odot u\right)  $ can be
presented as a matrix $%
\begin{array}
[c]{ccc}
& \Downarrow & \Downarrow\\
\Rightarrow & x & y\\
\Rightarrow & z & u
\end{array}
$ and for groups coincides with commutativity.

\begin{definition}
A ternary groupoid $(G,[\;])$ is \emph{medial} if it satisfies the identity
\[
\lbrack\lbrack x_{11}x_{12}x_{13}][x_{21}x_{22}x_{23}][x_{31}x_{32}%
x_{33}]]=[[x_{11}x_{21}x_{31}][x_{12}x_{22}x_{32}][x_{13}x_{23}x_{33}]]
\]
or%
\begin{equation}
m^{\left(  3\right)  }\circ\left(  m^{\left(  3\right)  }\times m^{\left(
3\right)  }\times m^{\left(  3\right)  }\right)  =m^{\left(  3\right)  }%
\circ\left(  m^{\left(  3\right)  }\times m^{\left(  3\right)  }\times
m^{\left(  3\right)  }\right)  \circ\sigma_{medial}, \label{smm}%
\end{equation}
where $\sigma_{medial}=\binom{123456789}{147258369}\in S_{9}.$
\end{definition}

It is not difficult to see that a semicommutative ternary semigroup is medial.

An element $x$ such that $[xxx]=x$ is called an \emph{idempotent}. A groupoid
in which all elements are idempotents is called an \emph{idempotent groupoid}.
A left (right, middle) identity is an idempotent.

\section*{Ternary groups and algebras}

\begin{definition}
A ternary semigroup $(G,[\;])$ is a \emph{ternary group} if for all $a,b,c\in
G$ there are $x,y,z\in G$ such that
\[
\lbrack xab]=[ayb]=[abz]=c.
\]
\end{definition}

One can prove \cite{pos}
%
that elements $x,y,z$ are uniquely determined. Moreover, according to the
suggestion of \cite{pos} one can prove (cf. \cite{dud/gla/gle})
%
that in the above definition, under the assumption of the associativity, it
suffices only to postulate the existence of a solution of $[ayb]=c$, or
equivalently, of $[xab]=[abz]=c$.

In a ternary group the equation $[xxz]=x$ has a unique solution which is
denoted by $z=\overline{x}$ and called \emph{skew element} (cf. \cite{dor3}),
or equivalently%
\[
m^{\left(  3\right)  }\circ\left(  \operatorname*{id}\times\operatorname*{id}%
\times\overline{\cdot}\right)  \circ D^{\left(  3\right)  }=\operatorname*{id}%
,
\]
where $D^{\left(  3\right)  }\left(  x\right)  =\left(  x,x,x\right)  $ is a
ternary diagonal map. As a consequence of results obtained in \cite{dor3} we have

\begin{theorem}
In any ternary group $(G,[\;])$ for all $x,y,z\in G$ the following relations
take place
\begin{align*}
\lbrack xx\,\overline{x}]  &  =[x\,\overline{x}\,x]=[\,\overline{x}\,xx]=x,\\
\lbrack yx\,\overline{x}]  &  =[y\,\overline{x}\,x]=[x\,\overline
{x}\,y]=[\,\overline{x}\,xy]=y,\\
\overline{\lbrack xyz]}  &  =[\,\overline{z}\,\overline{y}\,\overline{x}],\\
\overline{\overline{x}}  &  =x
\end{align*}
\end{theorem}

Since in an idempotent ternary group $\overline{x}=x$ for all $x$, an
idempotent ternary group is semicommutative. From results obtained in
\cite{dud/gla/gle} (see also \cite{dud2}) for $n=3$ we obtain

\begin{theorem}
A ternary semigroup $(G,[\;])$ with a unary operation ${}^{-}:x\rightarrow
\overline{x}$ is a ternary group if and only if it satisfies identities
\[
\lbrack yx\,\overline{x}\,]=[x\,\overline{x}\,y]=y,
\]
or%
\begin{align*}
m^{\left(  3\right)  }\circ\left(  \operatorname*{id}\times\overline{\cdot
}\times\operatorname*{id}\right)  \circ\left(  D^{\left(  2\right)  }%
\times\operatorname*{id}\right)   &  =\Pr\nolimits_{2},\\
m^{\left(  3\right)  }\circ\left(  \operatorname*{id}\times\operatorname*{id}%
\times\overline{\cdot}\right)  \circ\left(  \operatorname*{id}\times
D^{\left(  2\right)  }\right)   &  =\Pr\nolimits_{1},
\end{align*}
where $D^{\left(  2\right)  }\left(  x\right)  =\left(  x,x\right)  $ and
$\Pr\nolimits_{1}\left(  x,y\right)  =x$, $\Pr\nolimits_{2}\left(  x,y\right)
=y$.
\end{theorem}

\begin{corollary}
A ternary semigroup $(G,[\;])$ is an idempotent ternary group if and only if
it satisfies identities
\[
\lbrack yxx]=[xxy]=y
\]
\end{corollary}

A ternary group with an identity is derived from a binary group.

\begin{remark}
The set $A_{3}\subset S_{3}$ with ternary operation $[\;]$ defined as
composition of three permutations is an example of a ternary group which is
not derived from any group (all groups with three elements are commutative and
isomorphic to $\mathbb{Z}_{3}$).
\end{remark}

From results proved in \cite{dud2} follows

\begin{theorem}
A ternary group $(G, [\; ])$ satisfying the identity
\[
[xy\overline{x}]=y
\]
or
\[
[\overline{x}yx]=y
\]
is commutative.
\end{theorem}

\begin{theorem}
[Gluskin-Hossz\'{u}]For a ternary group $\left(  G,\left[  \;\right]  \right)
$ there exists a binary group $\left(  G,\circledast\right)  $, its
automorphism $\varphi$ and fixed element $b\in G$ such that
\begin{equation}
\left[  xyz\right]  =x\circledast\varphi\left(  y\right)  \circledast
\varphi^{2}\left(  z\right)  \circledast b.\label{gh}%
\end{equation}
\end{theorem}

\begin{proof}
Let $a\in G$ be fixed. Then the binary operation $x\circledast y=\left[
xay\right]  $ is associative, because
\[
(x\circledast y)\circledast z=[[xay]az]=[xa[yaz]]=x\circledast(y\circledast
z).
\]
An element $\overline{a}$ is its identity. $x^{-1}$ (in $(G,\circledast)$ is
$[\overline{a},\,\overline{x}\,\overline{a}]$. $\varphi(x)=\left[
\overline{a}xa\right]  $ is an automorphism of $(G,\circledast)$. The easy
calculation proves that the above formula holds for $b=[\overline
{a}\,\overline{a}\,\overline{a}\,]$. (see \cite{sok}).
\end{proof}

%
One can prove that the group $(G,\circledast)$ is unique up to isomorphism.
From the proof of Theorem 3 in \cite{gla/gle}
%
it follows that any medial ternary group satisfies the identity
\[
\overline{\lbrack xyz]}=[\overline{x}\,\overline{y}\,\overline{z}\,],
\]
which together with our previous results shows that in such groups we have
\[
\overline{\lbrack xyz]}=\overline{[xyz]}.
\]
But $\overline{\overline{x}}=x$. Hence, any medial ternary group is
semicommutative. Thus any retract of such group is a commutative group.
Moreover, for $\varphi$ from the proof of the above theorem we have
\[
\varphi(\varphi(x))=[\overline{a}\,[\overline{a}xa]a]=[\overline
{a}\,a\,[x\overline{a}\,a]]=x
\]

\begin{corollary}
Any medial ternary group $(G,[\;])$ has the form
\[
\lbrack xyz]=x\odot\varphi(y)\odot z\odot b,
\]
where $(G,\odot)$ is a commutative group, $\varphi$ its automorphism such that
$\varphi^{2}=\operatorname*{id}$ and $b\in G$ is fixed.
\end{corollary}

\begin{corollary}
A ternary group is medial if and only if it is semicommutative.
\end{corollary}

\begin{corollary}
A ternary group is semicommutative (medial) if and only if $[xay]=[yax]$ holds
for all $x,y\in G$ and some fixed $a\in G$.
\end{corollary}

\begin{corollary}
A commutative ternary group is $b$-derived from some commutative group.
\end{corollary}

Indeed, $\varphi(x)=[\overline{a}\,xa]=[xa\,\overline{a}]=x$.

\begin{theorem}
[Post]\label{theor-post}For any ternary group $\left(  G,\left[  \;\right]
\right)  $ there exists a binary group $\left(  G^{\ast},\circledast\right)  $
and $H\vartriangleleft G^{\ast}$, such that $G^{\ast}\diagup H\simeq
\mathbb{Z}_{2}$ and
\[
\left[  xyz\right]  =x\circledast y\circledast z
\]
for all $x,y,z\in G$.
\end{theorem}

\begin{proof}
Let $c$ be a fixed element in $G$ and let $G^{\ast}=G\times\mathbb{Z}_{2}$. In
$G^{\ast}$ we define binary operation $\circledast$ putting
\[
(x,0)\circledast(y,0) = ([xy\overline{c}], 1 )
\]
\[
(x,0)\circledast(y,1) = ([xyc], 0 )
\]
\[
(x,1)\circledast(y,0) = ([xcy], 0 )
\]
\[
(x,1)\circledast(y,1) = ([xcy], 1 ).
\]

It is not difficult to see that this operation is associative and
$(\overline{c},1)$ is its neutral element. The inverse element (in $G^{\ast}$)
has the form:
\[
(x,0)^{-1} = (\overline{x},0)
\]
\[
(x,1)^{-1} = ([\overline{c}\,\overline{x}\,\overline{c}], 1)
\]

Thus $G^{\ast}$ is a group such that $H=\{(x,1):x\in G\}\vartriangleleft
G^{\ast}$. Obviously the set $G$ can be identified with $G\times\{0\}$ and
\[
\lbrack xyz]=((x,0)\circledast(y,0))\circledast(z,0)=([xy\overline
{c}],1)\circledast(z,0)=
\]%
\[
([[xy\overline{c}]cz],0)=([xy[\overline{c}cz]],0)=([xyz],0)
\]
which completes the proof.
\end{proof}

\begin{proposition}
All retracts of a ternary group are isomorphic%
\[
ret_{a}\left(  G,\left[  \;\right]  \right)  \simeq ret_{b}\left(  G,\left[
\;\right]  \right)  .
\]
\end{proposition}

\begin{definition}
Autodistributivity in a ternary group is%
\[
\left[  \left[  xyz\right]  ab\right]  =\left[  \left[  xab\right]  \left[
yab\right]  \left[  zab\right]  \right]  .
\]
\end{definition}

Let us consider ternary algebras. Take 2 ternary operations $\left\{
\;,\;,\;\right\}  $ and $\left[  \;,\;,\;\right]  $, then distributivity is%
\[
\left\{  \left[  xyz\right]  ab\right\}  =\left[  \left\{  xab\right\}
\left\{  yab\right\}  \left\{  zab\right\}  \right]  ,
\]
and additivity is%
\[
\left[  \left\{  x+z\right\}  ab\right]  =\left[  xab\right]  +\left[
zab\right]  .
\]

\begin{definition}
Ternary algebra is a pair $\left(  A,m^{\left(  3\right)  }\right)  $, where
$A$ is a linear space and $m^{\left(  3\right)  }$ is a linear map
\[
m^{\left(  3\right)  }:A\otimes A\otimes A\rightarrow A
\]
called \textit{ternary multiplication} which is associative%
\[
m^{\left(  3\right)  }\circ\left(  m^{\left(  3\right)  }\otimes
\operatorname*{id}\otimes\operatorname*{id}\right)  =m^{\left(  3\right)
}\circ\left(  \operatorname*{id}\otimes m^{\left(  3\right)  }\otimes
\operatorname*{id}\right)  =m^{\left(  3\right)  }\circ\left(
\operatorname*{id}\otimes\operatorname*{id}\otimes m^{\left(  3\right)
}\right)  .
\]
\end{definition}

\section*{Ternary coalgebras}

Let $C$ is a linear space over a field $\mathbb{K}$.

\begin{definition}
Ternary comultiplication $\Delta^{\left(  3\right)  }$ is a linear map over a
fixed field $\mathbb{K}$%
\[
\Delta^{\left(  3\right)  }:C\rightarrow C\otimes C\otimes C.
\]
\end{definition}

For convenience we also use the short-cut Sweedler-type notations
\cite{sweedler}%
\begin{equation}
\Delta^{\left(  3\right)  }\left(  a\right)  =\sum_{i=1}^{n}a_{i}^{\prime
}\otimes a_{i}^{\prime\prime}\otimes a_{i}^{\prime\prime\prime}=a_{\left(
1\right)  }\otimes a_{\left(  2\right)  }\otimes a_{\left(  3\right)  }.
\label{sw}%
\end{equation}

Now we discuss various properties of $\Delta^{\left(  3\right)  }$ which are
in sense (dual) analog of the above ternary multiplication $m^{\left(
3\right)  }$.

First consider different possible types of ternary coassociativity.

\begin{enumerate}
\item \textit{Standard} ternary coassociativity%
\begin{equation}
(\Delta^{\left(  3\right)  }\otimes\operatorname*{id}\otimes\operatorname*{id}%
)\circ\Delta^{\left(  3\right)  }=(\operatorname*{id}\otimes\Delta^{\left(
3\right)  }\otimes\operatorname*{id})\circ\Delta^{\left(  3\right)
}=(\operatorname*{id}\otimes\operatorname*{id}\otimes\Delta^{\left(  3\right)
})\circ\Delta^{\left(  3\right)  }, \label{ast}%
\end{equation}

In the Sweedler notations%
\begin{align*}
&  \left(  a_{\left(  1\right)  }\right)  _{\left(  1\right)  }\otimes\left(
a_{\left(  1\right)  }\right)  _{\left(  2\right)  }\otimes\left(  a_{\left(
1\right)  }\right)  _{\left(  3\right)  }\otimes a_{\left(  2\right)  }\otimes
a_{\left(  3\right)  }=a_{\left(  1\right)  }\otimes\left(  a_{\left(
2\right)  }\right)  _{\left(  1\right)  }\otimes\left(  a_{\left(  2\right)
}\right)  _{\left(  2\right)  }\otimes\left(  a_{\left(  2\right)  }\right)
_{\left(  3\right)  }\otimes a_{\left(  3\right)  }\\
&  =a_{\left(  1\right)  }\otimes a_{\left(  2\right)  }\otimes\left(
a_{\left(  3\right)  }\right)  _{\left(  1\right)  }\otimes\left(  a_{\left(
3\right)  }\right)  _{\left(  2\right)  }\otimes\left(  a_{\left(  3\right)
}\right)  _{\left(  3\right)  }\equiv a_{\left(  1\right)  }\otimes a_{\left(
2\right)  }\otimes\otimes a_{\left(  3\right)  }\otimes a_{\left(  4\right)
}\otimes a_{\left(  5\right)  }.
\end{align*}

\item \textit{Nonstandard} ternary $\Sigma$-coassociativity (Gluskin-type ---
positional operatives)%
\[
(\Delta^{\left(  3\right)  }\otimes\operatorname*{id}\otimes\operatorname*{id}%
)\circ\Delta^{\left(  3\right)  }=(\operatorname*{id}\otimes\left(
\sigma\circ\Delta^{\left(  3\right)  }\right)  \otimes\operatorname*{id}%
)\circ\Delta^{\left(  3\right)  },
\]
where
\[
\sigma\circ\Delta^{\left(  3\right)  }\left(  a\right)  =\Delta_{\sigma
}^{\left(  3\right)  }\left(  a\right)  =a_{\left(  \sigma\left(  1\right)
\right)  }\otimes a_{\left(  \sigma\left(  2\right)  \right)  }\otimes
a_{\left(  \sigma\left(  3\right)  \right)  }%
\]
and $\sigma\in\Sigma\subset S_{3}.$

\item \textit{Permutational} ternary coassociativity%
\[
(\Delta^{\left(  3\right)  }\otimes\operatorname*{id}\otimes\operatorname*{id}%
)\circ\Delta^{\left(  3\right)  }=\pi\circ(\operatorname*{id}\otimes
\Delta^{\left(  3\right)  }\otimes\operatorname*{id})\circ\Delta^{\left(
3\right)  },
\]
where $\pi\in\Pi\subset S_{5}$.
\end{enumerate}

\textit{Ternary} \textit{comediality} is
\[
\left(  \Delta^{\left(  3\right)  }\otimes\Delta^{\left(  3\right)  }%
\otimes\Delta^{\left(  3\right)  }\right)  \circ\Delta^{\left(  3\right)
}=\sigma_{medial}\circ\left(  \Delta^{\left(  3\right)  }\otimes
\Delta^{\left(  3\right)  }\otimes\Delta^{\left(  3\right)  }\right)
\circ\Delta^{\left(  3\right)  },
\]
where $\sigma_{medial}$ is defined in (\ref{smm}).

\textit{Ternary counit} is defined as a map $\varepsilon^{\left(  3\right)
}:C\rightarrow\mathbb{K}$. In general, $\varepsilon^{\left(  3\right)  }%
\neq\varepsilon^{\left(  2\right)  }$ satisfying one of the conditions below.
If $\Delta^{\left(  3\right)  }$ is derived, then maybe $\varepsilon^{\left(
3\right)  }=\varepsilon^{\left(  2\right)  }$, but another counits may exist.

\begin{example}
Define $\left[  xyz\right]  =\left(  x+y+z\right)  |_{\operatorname{mod}\,2}$
for $x, y, z \in\mathbb{Z}_{2}$. It is seen that here there are 2 ternary
counits $\varepsilon^{\left(  3\right)  }=0,1$.
\end{example}

There are 3 types of ternary counits:

\begin{enumerate}
\item Standard (strong) ternary counit
\begin{equation}
(\varepsilon^{\left(  3\right)  }\otimes\varepsilon^{\left(  3\right)
}\otimes\operatorname*{id})\circ\Delta^{\left(  3\right)  }=(\varepsilon
^{\left(  3\right)  }\otimes\operatorname*{id}\otimes\varepsilon^{\left(
3\right)  })\circ\Delta^{\left(  3\right)  }=(\operatorname*{id}%
\otimes\varepsilon^{\left(  3\right)  }\otimes\varepsilon^{\left(  3\right)
})\circ\Delta^{\left(  3\right)  }=\operatorname*{id}, \label{e1}%
\end{equation}

\item Two sequensional (polyadic) counits $\varepsilon_{1}^{\left(  3\right)
}$ and $\varepsilon_{2}^{\left(  3\right)  }$%
\begin{equation}
(\varepsilon_{1}^{\left(  3\right)  }\otimes\varepsilon_{2}^{\left(  3\right)
}\otimes\operatorname*{id})\circ\Delta=(\varepsilon_{1}^{\left(  3\right)
}\otimes\operatorname*{id}\otimes\varepsilon_{2}^{\left(  3\right)  }%
)\circ\Delta=(\operatorname*{id}\otimes\varepsilon_{1}^{\left(  3\right)
}\otimes\varepsilon_{2}^{\left(  3\right)  })\circ\Delta=\operatorname*{id},
\label{e2}%
\end{equation}

\item Four long ternary counits $\varepsilon_{1}^{\left(  3\right)  }%
$--$\varepsilon_{4}^{\left(  3\right)  }$ satisfying%
\begin{equation}
\left(  \left(  \operatorname*{id}\otimes\varepsilon_{3}^{\left(  3\right)
}\otimes\varepsilon_{4}^{\left(  3\right)  }\right)  \circ\Delta^{\left(
3\right)  }\circ\left(  (\operatorname*{id}\otimes\varepsilon_{1}^{\left(
3\right)  }\otimes\varepsilon_{2}^{\left(  3\right)  })\circ\Delta^{\left(
3\right)  }\right)  \right)  =\operatorname*{id} \label{e3}%
\end{equation}
\end{enumerate}

By analogy with (\ref{xxx}) $\sigma$-\textit{cocommutativity} is defined
as\textbf{ }$\sigma\circ\Delta^{\left(  3\right)  }=\Delta^{\left(  3\right)
}$.

\begin{definition}
Ternary coalgebra is a pair $\left(  C,\Delta^{\left(  3\right)  }\right)  $,
where $C$ is a linear space and $\Delta^{\left(  3\right)  }$ is a ternary
comultiplication which is coassociative in one of the above senses.
\end{definition}

We will consider below only first standard type of associativity (\ref{ast}).

Let $\left(  A,m^{\left(  3\right)  }\right)  $ is a ternary algebra and
$\left(  C,\Delta^{\left(  3\right)  }\right)  $ is a ternary coalgebra and
$f,g,h\in\operatorname*{Hom}\nolimits_{\mathbb{K}}\left(  C,A\right)  .$

\begin{definition}
Ternary convolution product is%
\begin{equation}
\left[  f,g,h\right]  _{\ast}=m^{\left(  3\right)  }\circ\left(  f\otimes
g\otimes h\right)  \circ\Delta^{\left(  3\right)  } \label{tconv}%
\end{equation}
or $\left[  f,g,h\right]  _{\ast}\left(  a\right)  =\left[  f\left(
a_{\left(  1\right)  }\right)  g\left(  a_{\left(  2\right)  }\right)
h\left(  a_{\left(  3\right)  }\right)  \right]  $.
\end{definition}

\begin{definition}
Ternary coalgebra is called \textit{derived}, if there exists a binary (usual,
see e.g. \cite{abe,sweedler}) coalgebra $\Delta^{\left(  2\right)
}:C\rightarrow C\otimes C$ such that (cf. \ref{mder}))
\begin{equation}
\Delta_{der}^{\left(  3\right)  }=\left(  \operatorname*{id}\otimes
\Delta^{\left(  2\right)  }\right)  \otimes\Delta^{\left(  2\right)  }.
\label{dder}%
\end{equation}
\end{definition}

The derived ternary and $n$-ary coalgebras were considered e.g. in
\cite{san/mur} and \cite{bal/rag} respectively.

\section*{Ternary Hopf algebras}

\begin{definition}
Ternary bialgebra $B$ is triple $\left(  B,m^{\left(  3\right)  }%
,\Delta^{\left(  3\right)  }\right)  $ for which $\left(  B,m^{\left(
3\right)  }\right)  $ is a ternary algebra and $\left(  B,\Delta^{\left(
3\right)  }\right)  $ is a ternary coalgebra and%
\begin{equation}
\Delta^{\left(  3\right)  }\circ m^{\left(  3\right)  }=m^{\left(  3\right)
}\circ\Delta^{\left(  3\right)  } \label{md}%
\end{equation}
\end{definition}

One can distinguish four kinds of ternary bialgebrs with respect to a ''being
derived'' property'':

\begin{enumerate}
\item $\Delta$-derived ternary bialgebra%
\[
\Delta^{\left(  3\right)  }=\Delta_{der}^{\left(  3\right)  }=\left(
\operatorname*{id}\otimes\Delta^{\left(  2\right)  }\right)  \circ
\Delta^{\left(  2\right)  }%
\]

\item $m$-derived ternary bialgebra%
\[
m_{der}^{\left(  3\right)  }=m_{der}^{\left(  3\right)  }=m^{\left(  2\right)
}\circ\left(  m^{\left(  2\right)  }\otimes\operatorname*{id}\right)
\]

\item Derived ternary bialgebra is simultaneously $m$-derived and $\Delta
$-derived ternary bialgebra.

\item Full ternary bialgebra is not derived.
\end{enumerate}

Now we define possible types of ternary antipodes using analogy with binary coalgebras.

\begin{definition}
\textit{Skew ternary antipod} is%
\begin{align*}
m^{\left(  3\right)  }\circ(S_{skew}^{\left(  3\right)  }\otimes
\operatorname*{id}\otimes\operatorname*{id})\circ\Delta^{\left(  3\right)  }
&  =m^{\left(  3\right)  }\circ(\operatorname*{id}\otimes S_{skew}^{\left(
3\right)  }\otimes\operatorname*{id})\circ\Delta^{\left(  3\right)  }\\
&  =m^{\left(  3\right)  }\circ(\operatorname*{id}\otimes\operatorname*{id}%
\otimes S_{skew}^{\left(  3\right)  })\circ\Delta^{\left(  3\right)
}=\operatorname*{id}%
\end{align*}
or in terms of the ternary convolution product (\ref{tconv})%
\[
\left[  S_{skew}^{\left(  3\right)  },\operatorname*{id},\operatorname*{id}%
\right]  _{\ast}=\left[  \operatorname*{id},S_{skew}^{\left(  3\right)
},\operatorname*{id}\right]  _{\ast}=\left[  \operatorname*{id}%
,\operatorname*{id},S_{skew}^{\left(  3\right)  }\right]  _{\ast
}=\operatorname*{id}.
\]
\end{definition}

\begin{definition}
\textit{Strong ternary antipod} is%
\begin{align*}
\left(  m^{\left(  2\right)  }\otimes\operatorname*{id}\right)  \circ
(\operatorname*{id}\otimes S_{strong}^{\left(  3\right)  }\otimes
\operatorname*{id})\circ\Delta^{\left(  3\right)  }  &  =1\otimes
\operatorname*{id},\\
\left(  \operatorname*{id}\otimes m^{\left(  2\right)  }\right)
\circ(\operatorname*{id}\otimes\operatorname*{id}\otimes S_{strong}^{\left(
3\right)  })\circ\Delta^{\left(  3\right)  }  &  =\operatorname*{id}\otimes1,
\end{align*}
where $1$ is a unit of algebra.
\end{definition}

\begin{definition}
Ternary coalgebra is derived, if $\Delta^{\left(  3\right)  }$ is derived.
\end{definition}

\begin{lemma}
If in a ternary coalgebra $\left(  C,\Delta^{\left(  3\right)  }\right)  $
there exists a linear map $\varepsilon^{\left(  3\right)  }:C\rightarrow
\mathbb{K}$ satisfying%
\begin{equation}
\left(  \varepsilon^{\left(  3\right)  }\otimes\operatorname*{id}%
\otimes\varepsilon^{\left(  3\right)  }\right)  \circ\Delta^{\left(  3\right)
}=\operatorname*{id}, \label{ee}%
\end{equation}
then $\exists\Delta^{\left(  2\right)  }$ such that%
\[
\Delta^{\left(  3\right)  }=\Delta_{der}^{\left(  3\right)  }=\left(
\operatorname*{id}\otimes\Delta^{\left(  2\right)  }\right)  \otimes
\Delta^{\left(  2\right)  }%
\]
\end{lemma}

\begin{definition}
If in ternary coalgebra
\[
\Delta^{\left(  3\right)  }\circ S=\tau_{13}\circ\left(  S\otimes S\otimes
S\right)  \circ\Delta^{\left(  3\right)  },
\]
where $\tau_{13}=\binom{123}{321}$, then it is called \textit{skew-involutive}.
\end{definition}

\begin{definition}
Ternary Hopf algebra is a ternary bialgebra with a ternary antipod of the
corresponding type, i.e. $\left(  H,m^{\left(  3\right)  },e^{\left(
3\right)  },\Delta^{\left(  3\right)  },S^{\left(  3\right)  }\right)  $.
\end{definition}

\begin{remark}
There are 8 types of associative ternary Hopf algebras and 4 types of medial
Hopf algebras. Also it can happen that there are several ternary units
$e_{i}^{\left(  3\right)  }$ and several ternary counits $\varepsilon
_{i}^{\left(  3\right)  }$ (see (\ref{e1})--(\ref{e3})), which makes number of
possible ternary Hopf algebras enormous.
\end{remark}

\begin{theorem}
For any a ternary Hopf algebra there exists a binary Hopf algebra,
automorphism $\phi$ and a linear map $\lambda,$ such that%
\begin{equation}
\Delta^{\left(  3\right)  }=\left(  \operatorname*{id}\otimes\phi
\otimes\operatorname*{id}\right)  \circ\left(  \Delta^{\left(  2\right)
}\otimes\operatorname*{id}\right)  \label{d32}%
\end{equation}
\end{theorem}

\begin{proof}
The binary coproduct is $\Delta^{\left(  2\right)  }=\left(
\operatorname*{id}\otimes\lambda\otimes\operatorname*{id}\right)  \circ
\Delta^{\left(  3\right)  }$ and $\Delta^{\left(  3\right)  }=\left(
\operatorname*{id}\otimes\operatorname*{id}\otimes\operatorname*{id}%
\otimes\lambda\right)  \circ\left(  \Delta^{\left(  2\right)  }\otimes
\Delta^{\left(  2\right)  }\right)  \circ\Delta^{\left(  2\right)  }.$
\end{proof}

The co-analog of the Post Theorem \ref{theor-post} is

\begin{theorem}
For any ternary Hopf algebra $\left(  H,\Delta^{\left(  3\right)  }\right)  $
there exists a binary Hopf algebra $\left(  H^{\ast},\Delta^{\left(  2\right)
}\right)  $and $\Delta^{\left(  3\right)  }=\Delta_{der}^{\left(  3\right)
}|_{H}$, such that $H\diagup H^{\ast}\simeq k\left(  \mathbb{Z}_{2}\right)  $
and%
\begin{equation}
\left(  \operatorname*{id}\otimes\operatorname*{id}\otimes\operatorname*{id}%
\right)  \circ\Delta^{\left(  3\right)  }=\left(  \operatorname*{id}%
\otimes\Delta^{\left(  2\right)  }\right)  \circ\Delta^{\left(  2\right)  }.
\label{dpost}%
\end{equation}
\end{theorem}

\section*{Examples}

\begin{example}
Ternary dual pair $k\left(  G\right)  $ (push-forward) and $\mathcal{F}\left(
G\right)  $ (pull-back) which are related by $k^{\ast}\left(  G\right)
\cong\mathcal{F}\left(  G\right)  $. Here $k\left(  G\right)  =span\left(
G\right)  $ is a ternary group ($G$ has a ternary product $\left[  \;\right]
_{G}$ or $m_{G}^{\left(  3\right)  }$) algebra over a field $k$. If $u\in
k\left(  G\right)  $ ($u=u^{i}x_{i},x_{i}\in G$), then $\left[  uvw\right]
_{k}=u^{i}v^{j}w^{l}\left[  x_{i}x_{j}x_{l}\right]  _{G}$ is associative, and
so $\left(  k\left(  G\right)  ,\left[  \;\right]  _{k}\right)  $ becomes a
ternary algebra. Define a ternary coproduct $\Delta_{k}^{\left(  3\right)
}:k\left(  G\right)  \rightarrow k\left(  G\right)  \otimes k\left(  G\right)
\otimes k\left(  G\right)  $ by $\Delta_{k}^{\left(  3\right)  }\left(
u\right)  =u^{i}x_{i}\otimes x_{i}\otimes x_{i}$ (derive and associative),
then $\Delta_{k}^{\left(  3\right)  }\left(  \left[  uvw\right]  _{k}\right)
=\left[  \Delta_{k}^{\left(  3\right)  }\left(  u\right)  \Delta_{k}^{\left(
3\right)  }\left(  v\right)  \Delta_{k}^{\left(  3\right)  }\left(  w\right)
\right]  _{k}$, and $k\left(  G\right)  $ is a ternary bialgebra. If we define
a ternary antipod by $S_{k}^{\left(  3\right)  }=u^{i}\bar{x}_{i}$, where
$\bar{x}_{i}$ is a skew element of $x_{i}$, then $k\left(  G\right)  $ becomes
a ternary Hopf algebra. In the dual case of functions $\mathcal{F}\left(
G\right)  :\left\{  \varphi:G\rightarrow k\right\}  $ a ternary product
$\left[  \;\right]  _{\mathcal{F}}$ or $m_{\mathcal{F}}^{\left(  3\right)  }$
(derive and associative) acts on $\psi\left(  x,y,z\right)  $ as $\left(
m_{\mathcal{F}}^{\left(  3\right)  }\psi\right)  \left(  x\right)
=\psi\left(  x,x,x\right)  $, and so $\mathcal{F}\left(  G\right)  $ is a
ternary algebra. Let $\mathcal{F}\left(  G\right)  \otimes\mathcal{F}\left(
G\right)  \otimes\mathcal{F}\left(  G\right)  \cong\mathcal{F}\left(  G\times
G\times G\right)  $, then we define a ternary coproduct $\Delta_{\mathcal{F}%
}^{\left(  3\right)  }:\mathcal{F}\left(  G\right)  \rightarrow\mathcal{F}%
\left(  G\right)  \otimes\mathcal{F}\left(  G\right)  \otimes\mathcal{F}%
\left(  G\right)  $ as $\left(  \Delta_{\mathcal{F}}^{\left(  3\right)
}\varphi\right)  \left(  x,y,z\right)  =\varphi\left(  \left[  xyz\right]
_{\mathcal{F}}\right)  $, which is derive and associative. Thus we can obtain
$\Delta_{\mathcal{F}}^{\left(  3\right)  }\left(  \left[  \varphi_{1}%
\varphi_{2}\varphi_{3}\right]  _{\mathcal{F}}\right)  =\left[  \Delta
_{\mathcal{F}}^{\left(  3\right)  }\left(  \varphi_{1}\right)  \Delta
_{\mathcal{F}}^{\left(  3\right)  }\left(  \varphi_{2}\right)  \Delta
_{\mathcal{F}}^{\left(  3\right)  }\left(  \varphi_{3}\right)  \right]
_{\mathcal{F}}$, and therefore $\mathcal{F}\left(  G\right)  $ is a ternary
bialgebra. If we define a ternary antipod by $S_{\mathcal{F}}^{\left(
3\right)  }\left(  \varphi\right)  =\varphi\left(  \bar{x}\right)  $, where
$\bar{x}$ is a skew element of $x$, then $\mathcal{F}\left(  G\right)  $
becomes a ternary Hopf algebra.
\end{example}

\begin{example}
Matrix representation. Possible non-derived matrix representations of the
ternary product can be done only by four-rank tensors: twicely covariant and
twicely contravariant and allow only 2 possibilities $A_{jk}^{oi}B_{oo}%
^{jl}C_{il}^{ko}$ and $A_{ok}^{ij}B_{io}^{ol}C_{il}^{ko}$ (where $o$ is any index).
\end{example}

\textbf{Acknowledgments}. One of the authors (S.D.) would like to thank Jerzy
Lukierski for kind hospitality at the University of Wroc\l aw, where this work
was initiated and begun.

 \end{document}